\documentclass[preprint,12pt]{elsarticle}

\usepackage{amssymb}

\newtheorem{theorem}{Theorem}[section]
\newtheorem{lemma}[theorem]{Lemma}
\newtheorem{corollary}[theorem]{Corollary}

\newtheorem{question}[theorem]{Question}
\newtheorem{definition}[theorem]{Definition}

\journal{}

\begin{document}

\begin{frontmatter}



\title{Borel reducibility and H\" older($\alpha$) embeddability between Banach spaces}


\author{Longyun Ding\footnote{Research partially supported by the National Natural Science Foundation of China (Grant No. 10701044)}}

\ead{dinglongyun@gmail.com}
\address{School of Mathematical Sciences and LPMC, Nankai University, Tianjin 300071, PR China}

\begin{abstract}

We investigate Borel reducibility between equivalence relations
$E(X,p)=X^{\Bbb N}/\ell_p(X)$'s where $X$ is a separable Banach
space. We show that this reducibility is related to the so called
H\"older$(\alpha)$ embeddability between Banach spaces. By using the
notions of type and cotype of Banach spaces, we present many results
on reducibility and unreducibility between $E(L_r,p)$'s and
$E(c_0,p)$'s for $r,p\in[1,+\infty)$.

We also answer a problem presented by Kanovei in the affirmative by
showing that $C({\Bbb R}^+)/C_0({\Bbb R}^+)$ is Borel bireducible to
${\Bbb R}^{\Bbb N}/c_0$.

\end{abstract}

\begin{keyword} Borel reducibility \sep H\" older($\alpha$)
embeddability \sep type \sep cotype

\MSC 03E15 \sep 46B20 \sep 47H99




\end{keyword}

\end{frontmatter}


\section{Introduction}

Borel reducibility hierarchy of equivalence relations on Polish
spaces becomes the main focus of invariant descriptive set theory.
Some important equivalence relations, for example, $\Bbb R^\Bbb
N/\ell_p$ and $\Bbb R^\Bbb N/c_0$, were investigated by R. Dougherty
and G. Hjorth. They showed that, for $p,q\in[1,+\infty)$,
$$\Bbb R^\Bbb N/\ell_p\le_B\Bbb R^\Bbb N/\ell_q\iff p\le q,$$
while $\Bbb R^\Bbb N/c_0$ and $\Bbb R^\Bbb N/\ell_p$ are Borel
incomparable (see \cite{DH} and \cite{hjorth}).

Let $X$ be a topological linear space and $Y$ a Borel linear
subspace of $X$, then $X/Y$ is a natural example of Borel
equivalence relation. We are interested in the Borel reducibility
between this kind of equivalence relations. One of the motivation
for this paper is a problem asked by Kanovei that whether $C(\Bbb
R^+)/C_0(\Bbb R^+)\sim_B\Bbb R/c_0$.

In this paper, we generalize Doutherty and Hjorth's results by
considering equivalence relations $X^\Bbb N/\ell_p(X)$'s, which will
be denoted by $E(X,p)$ where $X$ is a separable Banach space and
$p\in[1,+\infty)$. We show that Borel reducibility between this kind
of equivalence relations is related to the existence of
H\"older($\alpha$) embeddings.

\begin{theorem}
Let $X,Y$ be two separable Banach spaces, $p,q\in[1,+\infty)$. If
there exists a H\"older$(\frac{p}{q})$ embedding from $X$ to
$\ell_q(Y)$, then we have $E(X,p)\le_B E(Y,q)$.
\end{theorem}

On the other hand, via introducing a similar notion of finitely
H\"older$(\alpha)$ embeddability, we prove the following theorem.

\begin{theorem}
Let $X,Y$ be two separable Banach spaces, $p,q\in[1,+\infty)$. If
$E(X,p)\le_B E(Y,q)$, then $X$ finitely H\"older$(\frac{p}{q})$
embeds into $\ell_q(Y)$.
\end{theorem}

For investigating the notion of finitely H\"older($\alpha$)
embaddability, we pay attention to the famous type-cotype theory in
Banach space theory.

\begin{theorem}
Let $X,U$ be two infinite dimensional Banach spaces, $\alpha>0$. If
$X$ finitely H\" older$(\alpha)$ embeds into $U$, then we have
\begin{enumerate}
\item[(1)] $\alpha\le 1$;
\item[(2)] $\frac{p(X)}{p(U)}\ge\alpha$;
\item[(3)] $p(U)>1\Rightarrow q(X)\le q(U)$.
\end{enumerate}
\end{theorem}

We apply this theorem to classical Banach spaces to show many
results on reducibility and unreducibility between $E(L_r,p)$'s and
$E(c_0,p)$'s. For instance, we show that

\begin{theorem}
For $r,s\in[1,2]$ and $p,q\in[1,+\infty)$, if $s\le q$, then
$$E(L_r,p)\le_B E(L_s,q)\iff p\le
q,\;\frac{r}{p}\ge\frac{s}{q}.$$
\end{theorem}

Equivalence relations $E(X,0)=X^\Bbb N/c_0(X)$'s are also
considered. We answer Kanovei's problem in the affirmative by
showing that

\begin{theorem}Let $X$ be a separable Banach space, $p\in[1,+\infty)$. Then
$E(X,0)\sim_B\Bbb R^\Bbb N/c_0$, and $E(X,p)$ is not Borel
comparable with $\Bbb R^\Bbb N/c_0$.
\end{theorem}

Since the proofs of several results of this paper are cited from
\cite{DH}, we shall assume that the reader has a copy of \cite{DH}
handy.

The paper is organized as follows. In section 2 we recall some
notions in descriptive set theory and functional analysis. In
section 3 we answer Kanovei's problem. In section 4 we introduce the
notion of finitely H\" older($\alpha$) embeddability and prove
Theorem 1.1 and 1.2. In section 5 we prove Theorem 1.3. In section 6
we apply theorems from earlier sections to classical Banach spaces,
and derive several interesting corollaries. Finally section 7
contains some further remarks.

\section{Basic notation}

For a set $I$, we denote by $|I|$ the cardinal of $I$.

A topological space is called a {\it Polish space} if it is
separable and completely metrizable. Let $X,Y$ be Polish spaces and
$E,F$ equivalence relations on $X,Y$ respectively. A {\it Borel
reduction} from $E$ to $F$  is a Borel function $\theta:X\to Y$ such
that
$$(x,y)\in E\iff(\theta(x),\theta(y))\in F$$ for all $x,y\in X$. We
say that $E$ is {\it Borel reducible} to $F$, denoted $E\le_B F$, if
there is a Borel reduction from $E$ to $F$. If $E\le_B F$ and
$F\le_B E$, we say that $E$ and $F$ are {\it Borel bireducible} and
denote $E\sim_B F$. Similarly, we say that $E$ is {\it strictly
Borel reducible} to $F$, denoted $E<_B F$, if $E\le_B F$ but not
$F\le_B E$. We refer to \cite{BK}, \cite{gaobook} and \cite{kanovei}
for background on Borel reducibility.

For two metric spaces $(M,d),(M',d')$ and $\alpha>0$. We say that
$M$  {\it H\"older$(\alpha)$ embeds} into $M'$ if there exist
$T:M\to M'$, $A>0$ such that
$$\frac{1}{A}d(u,v)^\alpha\le d'(T(u),T(v))\le
Ad(u,v)^\alpha$$ for $u,v\in M$.

As usual, we denote $L_p[0,1]$ by $L_p$ for $p\in[1,+\infty)$ for
the sake of brevity. For a Banach space $X$, we denote by
$\ell_p(X)$ the Banach space whose underlying space is
$$\left\{x\in X^\Bbb N:\sum_{n\in\Bbb N}\|x(n)\|_X^p<+\infty\right\},$$
with the norm
$$\|x\|_{X,p}=\left(\sum_{n\in\Bbb N}\|x(n)\|_X^p\right)^{\frac{1}{p}}.$$
For $m\in\Bbb N$, we also denote by $\ell_p^m(X)$ the finite
dimensional space $(X^m,\|\cdot\|_{X,p})$ where
$\|s\|_{X,p}=\left(\sum_{n=1}^m\|s(n)\|_X^p\right)^{\frac{1}{p}}$
for $s\in X^m$. Note that $\ell_p(\Bbb R)=\ell_p$ and $\ell_p^m(\Bbb
R)=\ell_p^m$.

The notions of type and cotype of Banach spaces are powerful tools
for investigating the local character of Banach spaces. An infinite
dimensional Banach space $X$ is said to have (Rademacher) {\it type}
$p$ for some $1\le p\le 2$, if there is a constant $C<+\infty$ such
that for every $n$ and any sequence $(u_j)_{j=1}^n$ in $X$,
$$\left(\frac{1}{2^n}\sum_{\epsilon\in\{-1,1\}^n}\left\|\sum_{j=1}^n\epsilon_ju_j\right\|^p\right)^{\frac{1}{p}}\le
C\left(\sum_{j=1}^n\|u_j\|^p\right)^{\frac{1}{p}};$$ $X$ is said to
have (Rademacher) {\it cotype} $q$ for some $q\ge 2$, if there is a
constant $D>0$ such that for every $n$ and any sequence
$(u_j)_{j=1}^n$ in $X$,
$$\left(\frac{1}{2^n}\sum_{\epsilon\in\{-1,1\}^n}\left\|\sum_{j=1}^n\epsilon_ju_j\right\|^q\right)^{\frac{1}{q}}\ge
D\left(\sum_{j=1}^n\|u_j\|^q\right)^{\frac{1}{q}}.$$ The supremum of
all types and the infimum of all cotypes of $X$ are denoted by
$p(X)$ and $q(X)$ respectively. For $p\in[1,+\infty)$, it is well
known that
$$p(\ell_p)=p(L_p)=\min\{p,2\},\quad q(\ell_p)=q(L_p)=\max\{p,2\}.$$
For any Banach space $X$, $r\in[1,+\infty)$, we have
$$p(\ell_r(X))=\min\{p(X),r\},\quad q(\ell_r(X))=\max\{q(X),r\}.$$
(We can prove these formulas from Kahane's inequality, or see the
remark in \cite{T-J}, pp.16) For more details on type and cotype, we
refer to \cite{LT}.

Let $X$ be an infinite dimensional Banach space, $p\in[1,+\infty)$,
we say that $X$ contains $\ell_p^n$'s uniformly if for $n\in\Bbb N$
and $\varepsilon>0$, there is a linear embedding $T_n:\ell_p^n\to X$
such that $\|T_n\|\cdot\|T_n^{-1}\|\le 1+\varepsilon$. The
Maurey-Pisier Theorem \cite{MP} shows that $X$ contains
$\ell_{p(X)}^n$'s and $\ell_{q(X)}^n$'s uniformly.

\section{Borel bireducibility with $\Bbb R^\Bbb N/c_0$}

Kanovei's problem is the following question \cite{gao},
Question~7.5. Another version of this problem is the Question~16.7.2
of \cite{kanovei}.

We denote the space of all positive real numbers by $\Bbb R^+$.

\begin{question}[Kanovei] Let $C(\Bbb R^+)$ be the space of
continuous functions on $\Bbb R^+$. We define an equivalence
relation $E_K$ by
$$f E_K g\iff\lim_{t\to+\infty}(f(t)-g(t))=0$$
for $f,g\in C(\Bbb R^+)$. Is $E_K\sim_B\Bbb R^\Bbb N/c_0$?
\end{question}

Before answering this question, we consider a class of equivalence
relations similar to $E_K$.

\begin{definition}
Let $(M_n,d_n),n\in\Bbb N$ be a sequence of separable complete
metric spaces. We define an equivalence relation $E((M_n)_{n\in\Bbb
N},0)$ on $\prod_{n\in\Bbb N}M_n$ by
$$(x,y)\in E((M_n)_{n\in\Bbb N},0)\iff\lim_{n\to\infty}d_n(x(n),y(n))=0$$
for $x,y\in\prod_{n\in\Bbb N}M_n$. If $(M_n,d_n)=(M,d)$ for
$n\in\Bbb N$, we write $E(M,0)=E((M_n)_{n\in\Bbb N},0)$ for the sake
of brevity.
\end{definition}

This notion was firstly introduced by I. Farah \cite{farah}, denoted
${\rm D}(\langle M_n,d_n\rangle)$. Many results on ${\rm D}(\langle
M_n,d_n\rangle)$ were given in \cite{farah}, especially the case
named $c_0$-equalities that all sets $M_n$ are finite. By this
definition, we have $E(\Bbb R,0)=\Bbb R^\Bbb N/c_0$. We can see that
$E_K$ is Borel reducible to $E(C[0,1],0)$ with the reducing map
$\theta_0:C(\Bbb R^+)\to C[0,1]^\Bbb N$ defined as
$$\theta_0(f)(n)(t)=f(t+n+1)$$
for $n\in\Bbb N$ and $t\in[0,1]$. To answer Kanovei's problem in the
affirmative, it will suffice to show that $E(C[0,1],0)\le_B\Bbb
R^\Bbb N/c_0$. In order to do so, we need a theorem of I. Aharoni.

\begin{theorem}[Aharoni \cite{aharoni}]\label{aharoni}
There is a constant $K>0$ such that for any separable metric space
$(M,d)$, there is a map $T:M\to c_0$ satisfying
$$d(u,v)\le\|T(u)-T(v)\|_{c_0}\le Kd(u,v)$$
for every $u,v\in M$.
\end{theorem}

We denote $I_n=\left\{\frac{k}{2^n}:k\in\Bbb N,0\le k\le
2^n\right\}$.

\begin{theorem}\label{c0}
Let $(M_n,d_n),n\in\Bbb N$ be a sequence of separable complete
metric spaces, then we have
\begin{enumerate}
\item[(i)] $E((M_n)_{n\in\Bbb N},0)\le_B\Bbb R^\Bbb N/c_0$;
\item[(ii)] $E((M_n)_{n\in\Bbb N},0)\le_B E([0,1],0)$;
\item[(iii)] $E((M_n)_{n\in\Bbb N},0)\le_B E((I_n)_{n\in\Bbb N},0)$.
\end{enumerate}
\end{theorem}

{\bf Proof.} (i) By Aharoni's theorem, there are $K>0$ and maps
$T_n:M_n\to c_0$ such that
$$d_n(u,v)\le\|T_n(u)-T_n(v)\|_{c_0}\le Kd_n(u,v)$$
for $u,v\in M_n$. Fix a bijection $\langle\cdot,\cdot\rangle:\Bbb
N^2\to\Bbb N$. We define $\theta:\prod_{n\in\Bbb N}M_n\to\Bbb R^\Bbb
N$ by
$$\theta(x)(\langle n,m\rangle)=T_n(x(n))(m)$$
for $x\in\prod_{n\in\Bbb N}M_n$ and $n,m\in\Bbb N$. It is easy to
see that $\theta$ is continuous.

Now we check that $\theta$ is a reduction.

For every $x,y\in\prod_{n\in\Bbb N}M_n$, if $(x,y)\in
E((M_n)_{n\in\Bbb N},0)$, then
$$\lim_{n\to\infty}d_n(x(n),y(n))\to 0.$$
So $\forall\varepsilon>0\exists N\forall
n>N(d_n(x(n),y(n))<\varepsilon)$. Since
$\|T_n(x(n))-T_n(y(n))\|_{c_0}\le Kd_n(x(n),y(n))<K\varepsilon$, we
have
$$\forall n>N\forall m(|T_n(x(n))(m)-T_n(y(n))(m)|<K\varepsilon).$$
For $n\le N$, since $T_n(x(n)),T_n(y(n))\in c_0$, we have
$$\lim_{m\to\infty}|T_n(x(n))(m)-T_n(y(n))(m)|=0.$$
Therefore, for all
but finitely many $(n,m)$'s, we have
$$|\theta(x)(\langle n,m\rangle)-\theta(y)(\langle
n,m\rangle)|=|T_n(x(n))(m)-T_n(y(n))(m)|<K\varepsilon.$$ Thus
$$\lim_{\langle n,m\rangle\to\infty}|\theta(x)(\langle n,m\rangle)-\theta(y)(\langle
n,m\rangle)|=0.$$ It follows that $\theta(x)-\theta(y)\in c_0$.

On the other hand, for every $x,y\in\prod_{n\in\Bbb N}M_n$, if
$\theta(x)-\theta(y)\in c_0$, then
$$\forall\varepsilon>0\exists N\forall n>N\forall m(|\theta(x)(\langle
n,m\rangle)-\theta(y)(\langle n,m\rangle)|<\varepsilon).$$
Therefore,
$$\begin{array}{ll}d_n(x(n),y(n))&\le\|T_n(x(n))-T_n(y(n))\|_{c_0}\cr
&=\sup\limits_{m\in\Bbb N}|T_n(x(n))(m)-T_n(y(n))(m)|\cr
&=\sup\limits_{m\in\Bbb N}|\theta(x)(\langle
n,m\rangle)-\theta(y)(\langle
n,m\rangle)|\le\varepsilon.\end{array}$$ It follows that $(x,y)\in
E((M_n)_{n\in\Bbb N},0)$.

(ii) Denote
$$d'_n(u,v)=\frac{d_n(u,v)}{1+d_n(u,v)}.$$
It is easy to see that
$$\lim_{n\to\infty}d_n(x(n),y(n))=0\iff\lim_{n\to\infty}d'_n(x(n),y(n))=0.$$

So we may assume that $d_n(u,v)\le 1$ for all $n\in\Bbb N$. By the
same arguments in (i), there are $K>0$, maps $T_n:M_n\to c_0$ such
that
$$d_n(u,v)\le\|T_n(u)-T_n(v)\|_{c_0}\le Kd_n(u,v)\le K$$
for $u,v\in M_n$ and a reduction $\theta:\prod_{n\in\Bbb
N}M_n\to\Bbb R^\Bbb N$ as $\theta(x)(\langle
n,m\rangle)=T_n(x(n))(m)$ for $x\in\prod_{n\in\Bbb N}M_n$ and
$n,m\in\Bbb N$.

There are real numbers $r_{n,m}\,(n,m\in\Bbb N)$ such that
$$T_n(u)(m)\in[r_{n,m},r_{n,m}+K]$$ for $u\in M_n$. Now we can define a
reducing map $\theta':\prod_{n\in\Bbb N}M_n\to{[0,1]}^\Bbb N$ by
$$\theta'(x)(\langle n,m\rangle)=\frac{T_n(x(n))(m)-r_{n,m}}{K}$$
for $x\in\prod_{n\in\Bbb N}M_n$ and $n,m\in\Bbb N$.

(iii) For $x\in[0,1]^\Bbb N$ and $n\in\Bbb N$, denote
$$\vartheta(x)(n)=\frac{[x(n)2^n]}{2^n}.$$
Clearly $\vartheta$ is a Borel reducing map from $[0,1]^\Bbb
N\to\prod_{n\in\Bbb N}I_n$. \hfill$\Box$

\begin{corollary}
$E_K\sim_B\Bbb R^\Bbb N/c_0\sim_B E([0,1],0)\sim_B E((I_n)_{n\in\Bbb
N},0)$.
\end{corollary}

{\bf Proof.} $\Bbb R^\Bbb N/c_0\le_B E_K$ is trivial. The remaining
parts of the corollary follow from Theorem~\ref{c0}.  \hfill$\Box$

\section{Borel reducibility between $E(X,p)$'s}

R. Dougherty and G. Hjorth proved in \cite{DH} that, for $1\le
p<q<+\infty$,
$$\Bbb R^\Bbb N/\ell^p<_B\Bbb R^\Bbb N/\ell^q.$$
In \cite{hjorth}, it is also shown by G. Hjorth that, for
$p\in[1,+\infty)$, $\Bbb R^\Bbb N/\ell^p$ and $\Bbb R^\Bbb N/c_0$
are $\le_B$ incomparable.

In the same spirit as of $\Bbb R^\Bbb N/\ell^p$ and
$E((M_n)_{n\in\Bbb N},0)$, we introduce the following equivalence
relations.

\begin{definition}
Let $(M_n,d_n),n\in\Bbb N$ be a sequence of separable complete
metric spaces. For $p\in[1,+\infty)$, we define an equivalence
relation $E((M_n)_{n\in\Bbb N},p)$ on $\prod_{n\in\Bbb N}M_n$ by
$$(x,y)\in E((M_n)_{n\in\Bbb N},p)\iff\sum_{n\in\Bbb N}d_n(x(n),y(n))^p<+\infty$$
for $x,y\in\prod_{n\in\Bbb N}M_n$. If $(M_n,d_n)=(M,d)$ for every
$n\in\Bbb N$, we write $E(M,p)=E((M_n)_{n\in\Bbb N},p)$ for the sake
of brevity.
\end{definition}

In this section, we mainly consider Borel reducibility between
$E(X,p)$'s where $X$'s are Banach spaces. It is straightforward to
check that $$E(X,p)=X^\Bbb N/\ell_p(X).$$

The following theorem gives a sufficient condition for
$E((M_n)_{n\in\Bbb N},p)\le_B E(X,q)$. This condition is quite
complicated though we need only some special cases.

\begin{theorem}\label{reduce}
Let $X$ be a separable Banach space, $p,q\in[1,+\infty)$, and let
$(M_n,d_n),n\in\Bbb N$ be a sequence of separable complete metric
spaces. If there are $A,C,D>0$, a sequence of Borel maps
$T_n:M_n\to\ell_q(X)$ and two sequences of non-negative real numbers
$\varepsilon_n,\delta_n,n\in\Bbb N$ such that
\begin{enumerate}
\item[(1)] $\sum_{n\in\Bbb N}\varepsilon_n^p<+\infty,\sum_{n\in\Bbb N}\delta_n^q<+\infty$;
\item[(2)]
$d_n(u,v)<\varepsilon_n\Rightarrow\|T_n(u)-T_n(v)\|_{X,q}<\delta_n$;
\item[(3)] $d_n(u,v)>C\Rightarrow\|T_n(u)-T_n(v)\|_{X,q}>D$;
\item[(4)] $\varepsilon_n\le d_n(u,v)\le C\Rightarrow$
$$\frac{1}{A}d_n(u,v)^{\frac{p}{q}}\le\|T_n(u)-T_n(v)\|_{X,q}\le
Ad_n(u,v)^{\frac{p}{q}}.$$
\end{enumerate}
Then we have $$E((M_n)_{n\in\Bbb N},p)\le_B E(X,q).$$
\end{theorem}

Before proving Theorem~\ref{reduce}, we present several easy
corollaries of it.

\begin{corollary}\label{reduce1}
Let $X$ be a separable Banach space and $(M,d)$ a separable complete
metric space, $p,q\in[1,+\infty)$. If $M$ H\"older$(\frac{p}{q})$
embeds into $\ell_q(X)$, then we have $E(M,p)\le_B E(X,q)$.
\end{corollary}

\begin{corollary}
Let $X$ and $Y$ be two Banach spaces, $p,q\in[1,+\infty)$. If there
exist $A,c,d>0$ and a Borel map $T:X\to\ell_q(Y)$ satisfying that
\begin{enumerate}
\item[(1)] $\|u-v\|_X<c\Rightarrow\|T(u)-T(v)\|_{Y,q}<d$;
\item[(2)] $\|u-v\|_X\ge c\Rightarrow$
$$\frac{1}{A}\|u-v\|_X^{\frac{p}{q}}\le\|T(u)-T(v)\|_{Y,q}\le
A\|u-v\|_X^{\frac{p}{q}}.$$
\end{enumerate}
Then we have $E(X,p)\le_B E(Y,q)$.
\end{corollary}

{\bf Proof.} Denote
$\varepsilon_n=2^{-n}c,\delta_n=(2^{\frac{p}{q}})^{-n}d$, define
$T_n:X\to\ell_q(Y)$ by
$$T_n(u)=(2^{\frac{p}{q}})^{-n}T(2^nu).$$
Then the result follows from Theorem~\ref{reduce}. \hfill$\Box$

{\bf Proof of Theorem~\ref{reduce}.} Fix a bijection
$\langle\cdot,\cdot\rangle:\Bbb N^2\to\Bbb N$. We define
$\theta:\prod_{n\in\Bbb N}M_n\to X^\Bbb N$ by
$$\theta(x)(\langle n,m\rangle)=T_n(x(n))(m)$$
for $x\in\prod_{n\in\Bbb N}M_n$ and $n,m\in\Bbb N$. It is easy to
see that $\theta$ is Borel. By the definition we have
$$\begin{array}{ll}&\sum_{n,m\in\Bbb
N}\|\theta(x)(\langle n,m\rangle)-\theta(y)(\langle
n,m\rangle)\|_X^q\cr =&\sum_{n\in\Bbb N}\sum_{m\in\Bbb
N}\|T_n(x(n))(m)-T_n(y(n))(m)\|_X^q\cr =&\sum_{n\in\Bbb
N}\|T_n(x(n))-T_n(y(n))\|_{X,q}^q\, .\end{array}$$

For $x,y\in\prod_{n\in\Bbb N}M_n$, we split $\Bbb N$ into three sets
$$I_1=\{n\in\Bbb N:d_n(x(n),y(n))<\varepsilon_n\},$$
$$I_2=\{n\in\Bbb N:d_n(x(n),y(n))>C\},$$
$$I_3=\{n\in\Bbb N:\varepsilon_n\le d_n(x(n),y(n))\le C\}.$$
From (2) we have $$\sum_{n\in I_1}d_n(x(n),y(n))^p<\sum_{n\in
I_1}\varepsilon_n^p\le\sum_{n\in\Bbb N}\varepsilon_n^p<+\infty,$$
$$\sum_{n\in I_1}\|T_n(x(n))-T_n(y(n))\|_{X,q}^q<\sum_{n\in
I_1}\delta_n^q\le\sum_{n\in\Bbb N}\delta_n^q<+\infty;$$ from (3) we
have
$$\sum_{n\in I_2}d_n(x(n),y(n))^p>C^p|I_2|,$$
$$\sum_{n\in I_2}\|T_n(x(n))-T_n(y(n))\|_{X,q}^q>D^q|I_2|;$$
and from (4) we have
$$\frac{1}{A^q}\sum_{n\in I_3}d_n(x(n),y(n))^p\le\sum_{n\in
I_3}\|T_n(x(n))-T_n(y(n))\|_{X,q}^q\le A^q\sum_{n\in
I_3}d_n(x(n),y(n))^p.$$ Therefore,
$$\begin{array}{ll}&(x,y)\in E((M_n)_{n\in\Bbb N},p)\cr\iff &\sum_{n\in\Bbb
N}d_n(x(n),y(n))^p<+\infty\cr\iff &|I_2|<\infty,\,\sum_{n\in
I_3}d_n(x(n),y(n))^p<+\infty\cr\iff &|I_2|<\infty,\,\sum_{n\in
I_3}\|T_n(x(n))-T_n(y(n))\|_{X,q}^q<+\infty\cr\iff &\sum_{n\in\Bbb
N}\|T_n(x(n))-T_n(y(n))\|_{X,q}^q<+\infty\cr\iff &\sum_{n,m\in\Bbb
N}\|\theta(x)(\langle n,m\rangle)-\theta(y)(\langle
n,m\rangle)\|_X^q<+\infty\cr\iff &\theta(x)-\theta(y)\in\ell^q(X).$$
\end{array}$$
It follows that $E((M_n)_{n\in\Bbb N},p)\le_B E(X,q)$. \hfill$\Box$

For Borel reducibility between $E(X,p)$'s, we aim to present a
necessary condition which will be named finitely
H\"older$(\frac{p}{q})$ embeddability. Now we focus on the
equivalence relations $E((Z_n)_{n\in\Bbb N},p)$ where $Z_n,n\in\Bbb
N$ are a sequence of finite metric spaces.

The following lemma is due to R. Dougherty and G. Hjorth.

\begin{lemma}\label{DH}
Let $Y$ be a separable Banach space, $p\in\{0\}\cup[1,+\infty)$,
$q\in[1,+\infty)$, and let $(Z_n,d_n),n\in\Bbb N$ be a sequence of
finite metric space. Assume that $E((Z_n)_{n\in\Bbb N},p)\le_B
E(Y,q)$. Then there exist strictly increasing sequences of natural
numbers $(b_j)_{j\in\Bbb N},(l_j)_{j\in\Bbb N}$ and $T_j:Z_{b_j}\to
Y^{l_{j+1}-l_j}$ such that, for $x,y\in\prod_{j\in\Bbb N}Z_{b_j}$,
we have
$$(x,y)\in E((Z_{b_j})_{j\in\Bbb N},p)\iff\sum_{j\in\Bbb
N}\|T_j(x(j))-T_j(y(j))\|_{Y,q}^q<+\infty.$$
\end{lemma}

{\bf Proof.} The proof is, almost word for word, a copy of the proof
of \cite{DH}, Theorem 2.2, Claim (i)-(iii).  \hfill$\Box$

Let $X,Y$ be two separable Banach spaces, $p,q\in[1,+\infty)$.
Assume that $E(X,p)\le_B E(Y,q)$.

Fix a sequence of finite subsets $F_n\subseteq X,n\in\Bbb N$ such
that
$$\{0\}\subseteq F_0\subseteq F_1\subseteq\cdots\subseteq F_n\subseteq\cdots$$
and $\bigcup_{n\in\Bbb N}F_n$ is dense in $X$. For every $n\in\Bbb
N$, we denote
$$Z_n=\left\{u+\frac{i}{2^n}(v-u):u,v\in F_n,0\le i\le
2^n\right\}.$$ Then $F_n\subseteq Z_n$.

Since $Z_n\subseteq X$ is a sequence of finite metric spaces, we can
find $(b_j)_{j\in\Bbb N},(l_j)_{j\in\Bbb N}$ and $T_j:Z_{b_j}\to
Y^{l_{j+1}-l_j}$, as in Lemma~\ref{DH}. Then we have the following
lemma.

\begin{lemma} \label{fholder}
There exists an $m\in\Bbb N$ such that $\forall k\exists N\forall
j>N$, for $u,v\in F_{b_j}$, if $\frac{1}{k}\le\|u-v\|_X\le 1$, then
we have
$$\frac{1}{2^m}\|u-v\|_X^{\frac{p}{q}}\le\|T_j(u)-T_j(v)\|_{Y,q}\le
2^m\|u-v\|_X^{\frac{p}{q}}.$$
\end{lemma}

{\bf Proof.} Assume for contradiction that, for every $m$, $\exists
k_m\exists^{\infty}j\exists u_j,v_j\in F_{b_j}$ such that
$\frac{1}{k_m}\le\|u_j-v_j\|_X\le 1$ but either
$$\frac{1}{2^m}\|u_j-v_j\|_X^{\frac{p}{q}}>\|T_j(u_j)-T_j(v_j)\|_{Y,q}$$ or
$$\|T_j(u_j)-T_j(v_j)\|_{Y,q}>2^m\|u_j-v_j\|_X^{\frac{p}{q}}.$$

We define two subsets $I_1,I_2\subseteq\Bbb N$. For $m\in\Bbb N$, we
put $m\in I_1$ iff $\exists k_m\exists^{\infty}j\exists u_j,v_j\in
F_{b_j}$ satisfying that $\frac{1}{k_m}\le\|u_j-v_j\|_X\le 1$ and
$$\frac{1}{2^m}\|u_j-v_j\|_X^{\frac{p}{q}}>\|T_j(u_j)-T_j(v_j)\|_{Y,q};$$
and $m\in I_2$ iff $\exists k_m\exists^{\infty}j\exists u_j,v_j\in
F_{b_j}$ satisfying that $\frac{1}{k_m}\le\|u_j-v_j\|_X\le 1$ and
$$\|T_j(u_j)-T_j(v_j)\|_{Y,q}>2^m\|u_j-v_j\|_X^{\frac{p}{q}}.$$

From the assumption, we can see that $I_1\cup I_2=\Bbb N$. Now we
consider the following two cases.

{\sl Case 1.} $|I_1|=\infty$. Select a finite set
$J_1^m\subseteq\Bbb N$ for every $m\in I_1$ satisfying that
\begin{enumerate}
\item[(i)] $|J_1^m|\le k_m^p$;
\item[(ii)] $1\le\sum_{j\in J_1^m}\|u_j-v_j\|_X^p\le 2$;
\item[(iii)] if $m_1<m_2$, then $\max J_1^{m_1}<\min J_1^{m_2}$.
\end{enumerate}

Now we define $x,y\in\prod_{j\in\Bbb N}Z_{b_j}$ by
$$\left\{\begin{array}{ll}x(j)=u_j,y(j)=v_j,& j\in J_1^m,m\in I_1,\cr
x(j)=y(j)=0,& \mbox{otherwise.}\end{array}\right.$$ Then we have
$$\sum_{j\in\Bbb N}\|x(j)-y(j)\|_X^p=\sum_{m\in I_1}\sum_{j\in
J_1^m}\|u_j-v_j\|_X^p\ge\sum_{m\in I_1}1=+\infty,$$ so $(x,y)\notin
E((Z_{b_j})_{j\in\Bbb N},p)$. One the other hand, we have
$$\begin{array}{ll}\sum_{j\in\Bbb N}\|T_j(x(j))-T_j(y(j))\|_{Y,q}^q&=\sum_{m\in I_1}\sum_{j\in
J_1^m}\|T_j(u_j)-T_j(v_j)\|_{Y,q}^q\cr &<\sum_{m\in I_1}\sum_{j\in
J_1^{m}}\frac{1}{2^{mq}}\|u_j-v_j\|_X^p\cr &\le 2\sum_{m\in
I_1}\left(\frac{1}{2^q}\right)^m\cr &<+\infty,\end{array}$$
contradicting Lemma~\ref{DH}!

{\sl Case 2.} $|I_2|=\infty$. Select a finite set
$J_2^m\subseteq\Bbb N$ for every $m\in I_2$ satisfying that
\begin{enumerate}
\item[(i)] $|J_2^m|\le k_m^p$;
\item[(ii)] $1\le\sum_{j\in J_2^m}\|u_j-v_j\|_X^p\le 2$;
\item[(iii)] if $m_1<m_2$, then $\max J_2^{m_1}<\min J_2^{m_2}$;
\item[(iv)] for $j\in J_2^m$, we have $m\le b_j$.
\end{enumerate}

For $m\in I_2,j\in J_2^m$, since $m\le b_j$, by the definition of
$Z_{b_j}$ we have
$$u_j^i\stackrel{\rm Def}{=}u_j+\frac{i}{2^m}(v_j-u_j)\in Z_{b_j}\quad(i=0,1,\cdots,2^m).$$
The triangle inequality gives
$$\sum_{1\le i\le
2^m}\|T_j(u_j^{i-1})-T_j(u_j^i)\|_{Y,q}\ge\|T_j(u_j)-T_j(v_j)\|_{Y,q},$$
thus there is an $i(j)$ such that
$$\|T_j(u_j^{i(j)-1})-T_j(u_j^{i(j)})\|_{Y,q}\ge\frac{1}{2^m}\|T_j(u_j)-T_j(v_j)\|_{Y,q}.$$

Now we define $x,y\in\prod_{j\in\Bbb N}Z_{b_j}$ by
$$\left\{\begin{array}{ll}x(j)=u_j^{i(j)-1},y(j)=u_j^{i(j)},& j\in J_2^m,m\in I_2,\cr
x(j)=y(j)=0,& \mbox{otherwise.}\end{array}\right.$$ Then we have
$$\begin{array}{ll}\sum_{j\in\Bbb N}\|x(j)-y(j)\|_X^p&=\sum_{m\in I_2}\sum_{j\in
J_2^m}\|u_j^{i(j)-1}-u_j^{i(j)}\|_X^p\cr &=\sum_{m\in I_2}\sum_{j\in
J_2^m}\frac{1}{2^{mp}}\|u_j-v_j\|_X^p\cr &\le 2\sum_{m\in
I_2}(\frac{1}{2^p})^m\cr &<+\infty,
\end{array}$$
so $(x,y)\in E((Z_{b_j})_{j\in\Bbb N},p)$. One the other hand, we
have
$$\begin{array}{ll}\sum_{j\in\Bbb N}\|T_j(x(j))-T_j(y(j))\|_{Y,q}^q&=\sum_{m\in I_2}\sum_{j\in
J_2^m}\|T_j(u_j^{i(j)-1})-T_j(u_j^{i(j)})\|_{Y,q}^q\cr &=\sum_{m\in
I_2}\sum_{j\in
J_2^m}\left(\frac{1}{2^m}\|T_j(u_j)-T_j(v_j)\|_{Y,q}\right)^q\cr
&>\sum_{m\in I_2}\sum_{j\in J_2^{m}}\|u_j-v_j\|_X^p\cr &\ge
\sum_{m\in I_2}1\cr &=+\infty,\end{array}$$ contradicting Lemma
\ref{DH} again! \hfill$\Box$

\begin{definition}
For two metric spaces $(M,d),(M',d')$ and $\alpha>0$. We say that
$M$  {\it finitely H\"older$(\alpha)$ embeds} into $M'$ if there
exists $A>0$ such that for every finite subset $F\subseteq M$, there
is $T_F:F\to M'$ satisfying
$$\frac{1}{A}d(u,v)^\alpha\le d'(T_F(u),T_F(v))\le
Ad(u,v)^\alpha$$ for $u,v\in F$. While $\alpha=1$, we also say $M$
finitely Lipschitz embeds into $M'$.
\end{definition}

\begin{theorem}\label{freduce}
Let $X,Y$ be two separable Banach spaces, $p,q\in[1,+\infty)$. The
the following conditions are equivalent:
\begin{enumerate}
\item[(1)] $X$ finitely H\"older$(\frac{p}{q})$ embeds into $\ell_q(Y)$.
\item[(2)] For any sequence of finite subsets $(F_n)_{n\in\Bbb N}$ of
$X$, we have $$E((F_n)_{n\in\Bbb N},p)\le_B E(Y,q).$$
\end{enumerate}
\end{theorem}

{\bf Proof.} (1)$\Rightarrow$(2). Since $X$ finitely
H\"older$(\frac{p}{q})$ embeds into $Y$, we can find $A>0$,
$T_n:F_n\to\ell_q(Y)$ such that
$$\frac{1}{A}\|u-v\|_X^{\frac{p}{q}}\le \|T_n(u)-T_n(v)\|_{Y,q}\le
A\|u-v\|_X^{\frac{p}{q}}$$ for $u,v\in F_n$. Then $E((F_n)_{n\in\Bbb
N},p)\le_B E(Y,q)$ follows from Theorem~\ref{reduce}.

(2)$\Rightarrow$(1). Fix a sequence of finite subsets $F_n\subseteq
X,n\in\Bbb N$ such that
$$\{0\}\subseteq F_0\subseteq F_1\subseteq\cdots\subseteq F_n\subseteq\cdots$$
and $\bigcup_{n\in\Bbb N}F_n$ is dense in $X$. Let $(b_j)_{j\in\Bbb
N},(l_j)_{j\in\Bbb N}$ and $T_j:F_{b_j}\to Y^{l_{j+1}-l_j}$ be from
Lemma~\ref{fholder}. For convenience, we identify
$(Y^{l_{j+1}-l_j},\|\cdot\|_{Y,q})$ with a subspace of $\ell_q(Y)$.
Then $T_j$ becomes a map $F_{b_j}\to\ell_q(Y)$.

Let us consider an arbitrary finite subset $F\subseteq X$. We can
find two integers $c,d>0$ such that $$\frac{1}{c}\le\|u-v\|_X\le
d,$$ or equivalently,
$$\frac{1}{cd}\le\left\|\frac{u}{d}-\frac{v}{d}\right\|_X\le 1$$ for
any distinct $u,v\in F$.

For every $u\in F$, since $\bigcup_{j\in\Bbb N}F_{b_j}$ is dense in
$X$, there exists an $R(u)\in\bigcup_{j\in\Bbb N}F_{b_j}$ such that
$$\left\|\frac{u}{d}-R(u)\right\|_X<\frac{1}{4cd}.$$ Then for any
distinct $u,v\in F$, we have
$$d\|R(u)-R(v)\|_X<d\left(\left\|\frac{u}{d}-\frac{v}{d}\right\|_X+\frac{1}{2cd}\right)=\|u-v\|_X+\frac{1}{2c}\le
2\|u-v\|_X,$$ and
$$d\|R(u)-R(v)\|_X>d\left(\left\|\frac{u}{d}-\frac{v}{d}\right\|_X-\frac{1}{2cd}\right)=\|u-v\|_X-\frac{1}{2c}\ge
\frac{1}{2}\|u-v\|_X.$$

From Lemma~\ref{fholder}, there exist $m\in\Bbb N$ and a
sufficiently large $i$ such that
\begin{enumerate}
\item[(i)] $R(u)\in F_{b_i}$ for every $u\in F$;
\item[(ii)] for $u,v\in F_{b_i}$, if $\frac{1}{cd}\le\|u-v\|_X\le 1$, then
$$\frac{1}{2^m}\|u-v\|_X^{\frac{p}{q}}\le\|T_i(u)-T_i(v)\|_{Y,q}\le
2^m\|u-v\|_X^{\frac{p}{q}}.$$
\end{enumerate}

We define $T_F:F\to\ell_q(Y)$ by
$$T_F(u)=d^{\frac{p}{q}}T_i(R(u))$$
for $u\in F$. Then for any distinct $u,v\in F$ we have
$$\begin{array}{ll}\|T_F(u)-T_F(u)\|_{Y,q}&=d^{\frac{p}{q}}\|T_i(R(u))-T_i(R(v))\|_{Y,q}\cr
&\le 2^m(d\|R(u)-R(v)\|_X)^{\frac{p}{q}}\cr
&<2^{m+\frac{p}{q}}\|u-v\|_X^{\frac{p}{q}},\end{array}$$ and
$$\begin{array}{ll}\|T_F(u)-T_F(u)\|_{Y,q}&=d^{\frac{p}{q}}\|T_i(R(u))-T_i(R(v))\|_{Y,q}\cr
&\ge 2^{-m}(d\|R(u)-R(v)\|_X)^{\frac{p}{q}}\cr
&>2^{-(m+\frac{p}{q})}\|u-v\|_X^{\frac{p}{q}}.\end{array}$$ Thus
$A=2^{m+\frac{p}{q}}$ witness that $X$ finitely
H\"older$(\frac{p}{q})$ embeds into $\ell_q(Y)$. \hfill$\Box$

\section{Finitely H\" older$(\alpha)$ embeddability between Banach spaces}

It is not surprising that finitely H\" older$(\alpha)$ embeddability
is related to ultraproducts of Banach spaces. An ultrafilter
$\mathfrak A$ on $\Bbb N$ is called free if it does not contain any
finite set. Let $U$ be a Banach space. Consider the space
$\ell_{\infty}(U)$ of all bounded sequences $x\in U^\Bbb N$ with the
norm $\|x\|=\sup_{n\in\Bbb N}\|x(n)\|_U$. Its subspace
$N=\{x:\lim_{\mathfrak A}\|x(n)\|_U=0\}$ is closed. The ultraproduct
$(U)_{\mathfrak A}$ is the quotient space $\ell_{\infty}(U)/N$ with
the norm $\|(x)_{\mathfrak A}\|_{\mathfrak A}=\lim_{\mathfrak
A}\|x(n)\|_U$. For more details on ultraproducts in Banach space
theory, see \cite{heinrich}.

\begin{theorem}
Let $X,U$ be two Banach spaces, $\alpha>0$, and let ${\mathfrak A}$
be a free ultrafilter on $\Bbb N$. Then $X$ finitely H\"
older$(\alpha)$ embeds into $U$ iff $X$ H\" older$(\alpha)$ embeds
into $(U)_{\mathfrak A}$.
\end{theorem}

{\bf Proof.} ($\Rightarrow$). Fix a sequence of finite subsets
$F_n\subseteq X,n\in\Bbb N$ such that
$$\{0\}\subseteq F_0\subseteq F_1\subseteq\cdots\subseteq F_n\subseteq\cdots$$
and $\bigcup_{n\in\Bbb N}F_n$ is dense in $X$.  There are $A>0$ and
$T_n:F_n\to U$ such that
$$\frac{1}{A}\|u-v\|_X^\alpha\le\|T_n(u)-T_n(v)\|_U\le
A\|u-v\|_X^\alpha$$ for $u,v\in F_n$. Let $u\in\bigcup_{n\in\Bbb
N}F_n$, $m=\min\{n:u\in F_n\}$, we define
$$T(u)=(0,\cdots,0,T_m(u),T_{m+1}(u),\cdots)_{\mathfrak A}.$$
By the definition of the norm on $(U)_{\mathfrak A}$, it is easy to
check that
$$\frac{1}{A}\|u-v\|_X^\alpha\le\|T(u)-T(v)\|_{\mathfrak A}\le
A\|u-v\|_X^\alpha$$ for $u,v\in\bigcup_{n\in\Bbb N}F_n$. Since
$\bigcup_{n\in\Bbb N}F_n$ is dense in $X$, we can extend $T$ onto
$X$.

($\Leftarrow$). Let $T:X\to(U)_{\mathfrak A}$ be a
H\"older$(\alpha)$ embedding with the constant $A>0$.

Fix a finite subset $F\subseteq X$. For $u,v\in F$, since
$$\frac{1}{A}\|u-v\|_X^\alpha\le\|T(u)-T(v)\|_{\mathfrak A}=\lim_{\mathfrak A}\|T(u)(n)-T(v)(n)\|_U\le
A\|u-v\|_X^\alpha,$$ we have
$$I_{u,v}\stackrel{Def}{=}\left\{n:\frac{1}{A+1}\|u-v\|_X^\alpha\le\|T(u)(n)-T(v)(n)\|_U\le
(A+1)\|u-v\|_X^\alpha\right\}\in\mathfrak A.$$ Now we fix an
$m\in\bigcap_{u,v\in F}I_{u,v}$. For $u\in F$, we define
$T_F(u)=T(u)(m)$ as desired. \hfill$\Box$

By using this theorem, we can transfer the problem of finitely
Lipschitz embeddability to the existence of Lipschitz embeddings.
The latter was deeply studied in geometric nonlinear functional
analysis. But this method does not work while $\alpha\ne 1$, because
there is no more known result on the existence of H\"older$(\alpha)$
embeddings. Most recently, we employed other powerful tools, i.e.,
metric type and metric cotype, to solve this problem.

\begin{lemma}\label{alpha}
Let $X,U$ be two Banach spaces, $\alpha>0$. If $X$ finitely H\"
older$(\alpha)$ embeds into $U$, then $\alpha\le 1$.
\end{lemma}

{\bf Proof.} Fix an $e\in X$ such that $\|e\|_X=1$. Denote
$F_n=\left\{\frac{i}{n}e:0\le i\le n\right\}$. There exist $A>0$ and
$T_n:F_n\to U$ such that
$$\frac{1}{A}=\frac{1}{A}\|e-0\|_X^\alpha\le\|T_n(e)-T_n(0)\|_U\le A\|e-0\|_X^\alpha,$$ and for $1\le i\le n$,
$$\frac{1}{A}\left\|\frac{i}{n}e-\frac{i-1}{n}e\right\|_X^\alpha\le\left\|T_n\left(\frac{i}{n}e\right)-T_n\left(\frac{i-1}{n}e\right)\right\|_U\le
A\left\|\frac{i}{n}e-\frac{i-1}{n}e\right\|_X^\alpha=\frac{A}{n^\alpha}.$$
The triangle inequality gives
$$\|T_n(e)-T_n(0)\|_U\le\sum_{i=1}^{n}\left\|T_n\left(\frac{i}{n}e\right)-T_n\left(\frac{i-1}{n}e\right)\right\|_U.$$
Thus $\frac{1}{A}\le n\cdot\frac{A}{n^\alpha}$, i.e.,
$$\frac{1}{A^2}\le n^{1-\alpha},$$
Letting $n\to\infty$, yields $\alpha\le 1$. \hfill$\Box$

The notion of metric type was introduced in \cite{BMW}. Let $(M,d)$
be a metric space. For a map $H:\{0,1\}^n\to M$ and
$s,s'\in\{0,1\}^n$, an unordered pair $\{H(s),H(s')\}$ is called an
{\it edge} if $s$ and $s'$ are different at exactly one coordinate;
and it is called a {\it diagonal} if $s$ and $s'$ are different at
all $n$ coordinates. Denote by $E$ the set of all edges and by $D$
the set of all diagonals. Clearly, $|E|=n2^{n-1},|D|=2^{n-1}$.

\begin{definition}[J. Bougain, V. Milman, H. Wolfson]
Let $p\ge 1$. A metric space $(M,d)$ has {\it metric type} $p$ if
there is a constant $C$, such that for every $n$ and any
$H:\{0,1\}^n\to M$ the following inequality holds:
$$\left(\sum_D d(H(s),H(s'))^2\right)^{\frac{1}{2}}\le
Cn^{\frac{1}{p}-\frac{1}{2}}\left(\sum_E
d(H(s),H(s'))^2\right)^{\frac{1}{2}}.$$ where the sums range over
all the diagonals and all the edges respectively.
\end{definition}

\begin{theorem}\label{mtype}
Let $X,U$ be two infinite dimensional Banach spaces, $\alpha>0$. If
$X$ finitely H\" older$(\alpha)$ embeds into $U$, then
$$\frac{p(X)}{p(U)}\ge\alpha.$$
\end{theorem}

{\bf Proof.} The Maurey-Pisier Theorem \cite{MP} says that $X$
contains $\ell_{p(X)}^n$'s uniformly. Thus for every $n\in\Bbb N$,
there is a linear operator $L_n:\ell_{p(X)}^n\to X$ such that for
all $s\in\ell_{p(X)}^n$ we have $\|s\|_{p(X)}\le\|L_n(s)\|_X\le
2\|s\|_{p(X)}$.

Note that $\{0,1\}^n\subseteq\ell_{p(X)}^n$, we denote
$F_n=L_n(\{0,1\}^n)$. There are $A>0$ and $T_n:F_n\to U$ such that
for $u,v\in F_n$ we have
$$\frac{1}{A}\|u-v\|_X^\alpha\le\|T_n(u)-T_n(v)\|_U\le
A\|u-v\|_X^\alpha.$$

Fix a $p<p(U)$. By \cite{BMW}, Corollary 5.10, $U$ has metric type
$p$. Denote $H=T_n\circ L_n$. Now we estimate lengths of edges and
diagonals for $H$ as follows.

Let $s,s'\in\{0,1\}^n$. If $s$ and $s'$ are different at exactly one
coordinate, then $\|s-s'\|_{p(X)}=1,\|L_n(s)-L_n(s')\|_X\le 2$. So
$$\|H(s)-H(s')\|_U=\|T_n(L_n(s))-T_n(L_n(s'))\|_U\le A2^\alpha.$$
On the other hand, if $s$ and $s'$ are different at all $n$
coordinates, then
$\|s-s'\|_{p(X)}=n^{\frac{1}{p(X)}},\|L_n(s)-L_n(s')\|_X\ge
n^{\frac{1}{p(X)}}$. So
$$\|H(s)-H(s')\|_U=\|T_n(L_n(s))-T_n(L_n(s'))\|_U\ge\frac{n^{\frac{\alpha}{p(X)}}}{A}.$$
Therefore,
$$\begin{array}{ll}2^{n-1}n^{\frac{2\alpha}{p(X)}}A^{-2}&\le\sum_D\|H(s)-H(s')\|_U^2\cr
&\le C^2n^{2(\frac{1}{p}-\frac{1}{2})}\sum_E\|H(s)-H(s')\|_U^2\cr
&\le
C^2n^{2(\frac{1}{p}-\frac{1}{2})}n2^{n-1}A^22^{2\alpha}.\end{array}$$
Thus $$n^{\frac{\alpha}{p(X)}-\frac{1}{p}}\le CA^22^\alpha.$$

By letting $n\to\infty$, we see that
$\frac{\alpha}{p(X)}\le\frac{1}{p}$ for every $p<p(U)$. It then
follows that $\frac{\alpha}{p(X)}\le\frac{1}{p(U)}$, i.e.,
$\frac{p(X)}{p(U)}\ge\alpha$. \hfill$\Box$

The notion of metric cotype introduced by M. Mendel and A. Naor
\cite{MN} is more complicated than metric type.

\begin{definition}[M. Mendel, A. Naor]\label{cotype} Let $q>0$. A metric space
$(M,d)$ has {\it metric cotype} $q$ if there is a constant $\Gamma$,
which satisfies that for every $n$, there exists an even integer
$m$, such that for every $H:{\Bbb Z}_m^n\to M$, the following
inequality holds:
$$\sum_{j=1}^n\Bbb E_s\left[d\left(H\left(s+\frac{m}{2}e_j\right),H(s)\right)^q\right]\le\Gamma^qm^q\Bbb E_{\epsilon,s}[d(H(s+\epsilon),H(s))^q],$$
where the expectations $\Bbb E_s$ and $\Bbb E_{\epsilon,s}$ above
are taken with respect to uniformly chosen $s\in{\Bbb Z}_m^n$ and
$\epsilon\in\{-1,0,1\}^n$, and
$e_j=(0,\cdots,0,\stackrel{j}{1},0,\cdots,0)$ for $j=1,\cdots,n$.
\end{definition}

M. Mendel and A. Naor showed that, for a Banach space, metric cotype
is coincide with cotype (see \cite{MN}, Theorem 1.2). They also
estimated, in Theorem 4.1 of \cite{MN}, the minimal $m$ in
Definition~\ref{cotype} for $K$-convex Banach spaces.

\begin{lemma}[M. Mendel, A. Naor]\label{MN}
Let $X$ be a $K$-convex Banach space with cotype $q$. Then there
exists constant $C>0$ such that for every $n$ and every integer
$m\ge Cn^{\frac{1}{q}}$ which is divisible by $4$, the inequality in
Definition~\ref{cotype} holds.
\end{lemma}

\begin{theorem}\label{mcotype}
Let $X,U$ be two infinite dimensional Banach spaces with $p(U)>1$.
For $\alpha>0$, if $X$ finitely H\" older$(\alpha)$ embeds into $U$,
then $$q(X)\le q(U).$$
\end{theorem}

{\bf Proof.} For every $n,m$, the map $\sigma_n:{\Bbb Z}_m^n\to\Bbb
C^n$ is defined by
$$\sigma_n(k_1,\cdots,k_n)=\left(\exp\left(\frac{2\pi
k_1}{m}i\right),\cdots,\exp\left(\frac{2\pi k_n}{m}i\right)\right)$$
for $k_1,\cdots,k_n\in{\Bbb Z}_m$. Note that $\ell_{q(X)}^n(\Bbb
C)\cong\ell_{q(X)}^n(\Bbb R^2)$, so it is $L$-isomorphic to
$\ell_{q(X)}^{2n}$ where $L>0$ is a constant independent to $n$. By
the Maurey-Pisier Theorem \cite{MP}, $X$ contains $\ell_{q(X)}^n$'s
uniformly. Thus we can find a linear operation
$R_n:\ell_{q(X)}^n(\Bbb C)\to X$ and $P,Q>0$ independent of $n$ such
that
$$P\|s\|_{\Bbb C,q(X)}\le\|R_n(s)\|_X\le Q\|s\|_{\Bbb C,q(X)}$$
for $s\in\ell_{q(X)}^n(\Bbb C)$.

We denote $F_n=R_n(\sigma_n({\Bbb Z}_m^n))$. There are $A>0$ and
$T_n:F_n\to U$ such that for $u,v\in F_n$ we have
$$\frac{1}{A}\|u-v\|_X^\alpha\le\|T_n(u)-T_n(v)\|_U\le
A\|u-v\|_X^\alpha.$$

Fix a $q>q(U)$. Denote $H=T_n\circ R_n\circ\sigma_n$. By the
Pisier's $K$-convexity theorem \cite{pisier}, $p(U)>1$ iff $U$ is
$K$-convex. Since $U$ has cotype $q$, from Lemma~\ref{MN}, there is
a constant $C>0$ and for every sufficiently large $n$, we can find a
suitable $m$ such that $m$ is divisible by $4$, $Cn^{\frac{1}{q}}\le
m\le(C+1)n^{\frac{1}{q}}$ and
$$\sum_{j=1}^n\Bbb E_s\left[\left\|H\left(s+\frac{m}{2}e_j\right)-H(s)\right\|_U^q\right]\le\Gamma^qm^q\Bbb E_{\epsilon,s}[\|H(s+\epsilon)-H(s)\|_U^q].$$

For $s\in{\Bbb Z}_m^n$ and $\epsilon\in\{-1,0,1\}^n$, we have
$$\begin{array}{ll}\|\sigma_n(s+\epsilon)-\sigma_n(s)\|_{\Bbb C,q(X)}
&=\left(\sum_{j=1}^n\left|\exp\left(\frac{2\pi(s_j+\epsilon_j)}{m}i\right)-\exp\left(\frac{2\pi
s_j}{m}i\right)\right|^{q(X)}\right)^{\frac{1}{q(X)}}\cr
&\le\left(\sum_{j=1}^n\left|\frac{2\pi\epsilon_j}{m}\right|^{q(X)}\right)^{\frac{1}{q(X)}}\le\frac{2\pi}{m}n^{\frac{1}{q(X)}},\end{array}$$
so
$$\begin{array}{ll}\|H(s+\epsilon)-H(s)\|_U&\le
A\|R_n(\sigma_n(s+\epsilon))-R_n(\sigma_n(s))\|_X^\alpha\cr &\le
AQ^\alpha\|\sigma_n(s+\epsilon)-\sigma_n(s)\|_{\Bbb
C,q(X)}^\alpha\cr &\le AQ^\alpha(2\pi)^\alpha
n^{\frac{\alpha}{q(X)}}m^{-\alpha}.\end{array}$$ Moreover, for
$j=1,2,\cdots,n$, we have
$$\left\|\sigma_n\left(s+\frac{m}{2}e_j\right)-\sigma_n(s)\right\|_{\Bbb C,q(X)}=\left|\exp\left(\frac{2\pi s_j}{m}i+\pi
i\right)-\exp\left(\frac{2\pi s_j}{m}i\right)\right|=2,$$ so
$$\begin{array}{ll}\left\|H\left(s+\frac{m}{2}e_j\right)-H(s)\right\|_U
&\ge\frac{1}{A}\left\|R_n\left(\sigma_n\left(s+\frac{m}{2}e_j\right)\right)-R_n(\sigma_n(s))\right\|_X^\alpha\cr
&\ge
A^{-1}P^\alpha\left\|\sigma_n\left(s+\frac{m}{2}e_j\right)-\sigma_n(s)\right\|_{\Bbb
C,q(X)}^\alpha\cr &\ge 2^\alpha A^{-1}P^\alpha.\end{array}$$
Therefore,
$$\begin{array}{ll}n(2^\alpha A^{-1}P^\alpha)^q&\le\sum_{j=1}^n\Bbb E_s\left[\left\|H\left(s+\frac{m}{2}e_j\right)-H(s)\right\|_U^q\right]\cr
&\le\Gamma^qm^q(AQ^\alpha(2\pi)^\alpha
n^{\frac{\alpha}{q(X)}}m^{-\alpha})^q,\end{array}$$ i.e.,
$$n^{\frac{1}{q}-\frac{\alpha}{q(X)}}m^{\alpha-1}\le
W\stackrel{Def}{=}\Gamma A^2P^{-\alpha}Q^\alpha\pi^\alpha.$$
Lemma~\ref{alpha} gives $\alpha\le 1$. Since
$m\le(C+1)n^{\frac{1}{q}}$ for sufficiently large $n$, we have
$m^{\alpha-1}\ge(C+1)^{\alpha-1}n^{\frac{\alpha-1}{q}}$. Thus
$$n^{\frac{1}{q}-\frac{1}{q(X)}}\le(W(C+1)^{-(\alpha-1)})^{\frac{1}{\alpha}}.$$

By letting $n\to\infty$, we see that $\frac{1}{q}\le\frac{1}{q(X)}$
for every $q>q(U)$. It then follows that $q(X)\le q(U)$.
\hfill$\Box$

\section{Applications to classical Banach spaces}

In this section, we compare equivalence relations $E(X,p)$'s, where
$p\in\{0\}\cup[1,+\infty)$ and $X$ is one of classical Banach
spaces, namely, $X=c_0,C[0,1]$ or $X=\ell_r,L_r$ for
$r\in[1,+\infty)$.

Firstly, we present all reducibility concerning the case $p=0$.

\begin{theorem}
Let $X$ be a separable Banach space, $p\in[1,+\infty)$. Then
$E(X,0)\sim_B\Bbb R^\Bbb N/c_0$, and $E(X,p)$ is not Borel
comparable with $\Bbb R^\Bbb N/c_0$.
\end{theorem}

{\bf Proof.} $\Bbb R^\Bbb N/c_0\le_B E(X,0)$ is trivial, and
$E(X,0)\le_B\Bbb R^\Bbb N/c_0$ follows from Theorem~\ref{c0}.

R. Dougherty and G. Hjorth showed that $\Bbb R^\Bbb
N/\ell_1\not\le_B\Bbb R^\Bbb N/c_0$ and $\Bbb R^\Bbb
N/\ell_1\le_B\Bbb R^\Bbb N/\ell_p$ (see \cite{hjorth}, Theorem 6.1,
and \cite{DH}, Theorem 1.1). By Corollary~\ref{reduce1}, we have
$\Bbb R^\Bbb N/\ell_p=E(\Bbb R,p)\le_B E(X,p)$. Therefore,
$E(X,p)\not\le_B\Bbb R^\Bbb N/c_0$.

Suppose $\Bbb R^\Bbb N/c_0\le E(X,p)$. Denote
$Z_n=I_n=\left\{\frac{k}{2^n}:k\in\Bbb N,0\le k\le 2^n\right\}$. By
Lemma~\ref{DH}, there exist strictly increasing sequences of natural
numbers $(b_j)_{j\in\Bbb N},(l_j)_{j\in\Bbb N}$ and $T_j:Z_{b_j}\to
X^{l_{j+1}-l_j}$ such that, for $x,y\in\prod_{j\in\Bbb N}Z_{b_j}$,
we have
$$(x,y)\in\Bbb R^\Bbb N/c_0\iff\sum_{j\in\Bbb
N}\|T_j(x(j))-T_j(y(j))\|_{X,p}^p<+\infty.$$ By the same trick of
the proof of \cite{DH}, Theorem 2.2, Claim (iv), we obtain a
contradiction. \hfill$\Box$

\begin{lemma}\label{fractal}
For any $\alpha\in(0,1]$, there is an $n\in\Bbb N$ such that $\Bbb
R$ H\"older$(\alpha)$ embeds into $\Bbb R^n$.
\end{lemma}

{\bf Proof.} It will suffice to prove for the case
$\frac{1}{2}<\alpha<1$.

Let $r=4^{-\alpha}$. Then $\frac{1}{4}<r<\frac{1}{2}$. Proposition
1.2 of \cite{DH} gives a H\"older$(\alpha)$ embedding
$K_r:[0,1]\to\Bbb R^2$. As remarked in \cite{DH}, we can extend
$K_r$ to all of $\Bbb R$ as follows.

First step, we extend $K_r$ to a H\"older$(\alpha)$ embedding
$K_r^1:[0,4]\to\Bbb R^2$ by
$$K_r^1(t)=r^{-1}K_r\left(\frac{t}{4}\right)$$ for $t\in[0,4]$. From the definition
of $K_r(t)$, note that
$$K_r(0)=(0,0),K_r\left(\frac{1}{4}\right)=(r,0),K_r\left(\frac{3}{4}\right)=(1-r,0),K_r(1)=(1,0),$$
we can see that $K_r^1\upharpoonright[0,1]=K_r$.

Second step, we extend $K_r^1$ to a H\"older$(\alpha)$ embedding
$K_r^2:[-12,4]\to\Bbb R^2$ by
$$K_r^2(t)=r^{-1}K_r^1\left(\frac{t+12}{4}\right)-(r^{-2}-r^{-1},0)$$ for
$t\in[-12,4]$. From the definition of $K_r(t)$, note that
$$K_r^1(0)=(0,0),K_r^1(3)=(r^{-1}-1,0),K_r^1(4)=(r^{-1},0),$$
we can see that $K_r^2\upharpoonright[0,4]=K_r^1$.

Repeating these steps, we can extend $K_r$ to a H\"older$(\alpha)$
embedding $K_r^\infty:\Bbb R\to\Bbb R^2$.

For $\frac{1}{2^k}<\alpha<1$, by repeatedly applying $K_r^\infty$
for some suitable $r$, we can find a H\"older$(\alpha)$ embedding
$\Bbb R\to\Bbb R^{2^k}$. \hfill$\Box$

\begin{theorem}\label{reducible}
For $r,s,p,q\in[1,+\infty)$, we have
\begin{enumerate}
\item[(1)] $E(\ell_r,p)\le_B E(L_r,p)$;
\item[(2)] $E(\ell_p,p)\sim_B\Bbb R^\Bbb N/\ell_p\le_B E(\ell_r,p)$;
\item[(3)] $E(\ell_2,p)\sim_B E(L_2,p)\le_B E(L_r,p)$;
\item[(4)] if $s\le r\le 2$, then $E(L_r,p)\le_B E(L_s,p)$;
\item[(5)] if $\frac{r}{p}=\frac{s}{q},p\le q$, then $E(\ell_r,p)\le_B
E(\ell_s,q)$ and $E(L_r,p)\le_B E(L_s,q)$.
\end{enumerate}
\end{theorem}

{\bf Proof.} Let $X\hookrightarrow Y$ stand for that $X$ Lipschitz
embeds into $Y$. From Corollary~\ref{reduce1}, clauses (1)-(4) are
given by following well known facts.
\begin{enumerate}
\item[(1)] $\ell_r\hookrightarrow L_r\hookrightarrow\ell_p(L_r)$.
\item[(2)] $\Bbb R^\Bbb N/\ell_p=E(\Bbb R,p)$ and $\Bbb
R\hookrightarrow X$ for any Banach space.
\item[(3)] $\ell_2\cong L_2$ and $L_2\hookrightarrow L_r$ (see \cite{BL}, pp. 189).
\item[(4)] if $s\le r\le 2$, then $L_r\hookrightarrow L_s$ (see \cite{BL}, Corollary 8.8).
\end{enumerate}

(5) Denote $\alpha=\frac{r}{s}=\frac{p}{q}\in(0,1]$. By Lemma
\ref{fractal}, there are $n\in\Bbb N$ and a H\"older$(\alpha)$
embedding $T:\Bbb R\to\Bbb R^n$. Hence, there exists $A>0$, such
that for $t_1,t_2\in\Bbb R$, we have
$$\frac{1}{A}|t_1-t_2|^\alpha\le\|T(t_1)-T(t_2)\|_2\le A|t_1-t_2|^\alpha.$$

Now we define a H\"older$(\alpha)$ embedding $\tilde T:L_r\to L_s$.
For $f\in L_r$ and $t\in(0,1]$, if
$t\in\left(\frac{k-1}{n},\frac{k}{n}\right]$ for some
$k\in\{1,2,\cdots,n\}$, let
$$\tilde T(f)(t)=n^{\frac{1}{s}}T(f(nt-k+1))(k).$$
For $f,g\in L_r$, denote $\tau=nt-k+1$, then
$$\int_\frac{k-1}{n}^\frac{k}{n}|\tilde T(f)(t)-\tilde T(g)(t)|^s
dt=\int_0^1|T(f(\tau))(k)-T(g(\tau))(k)|^s d\tau.$$ Therefore,
$$\begin{array}{ll}\int_0^1|\tilde T(f)(t)-\tilde T(g)(t)|^s dt
&=\sum_{k=1}^{n}\int_\frac{k-1}{n}^\frac{k}{n}|\tilde T(f)(t)-\tilde
T(g)(t)|^s dt\cr
&=\int_0^1\sum_{k=1}^{n}|T(f(\tau))(k)-T(g(\tau))(k)|^s d\tau\cr
&=\int_0^1\|T(f(\tau))-T(g(\tau))\|_s^s d\tau.\end{array}$$ Since
$$\frac{1}{n}\|u\|_2\le\|u\|_\infty\le\|u\|_s\le n\|u\|_\infty\le n\|u\|_2,$$ for $u\in\Bbb
R^n$, we have
$$\begin{array}{ll}\int_0^1|\tilde T(f)(t)-\tilde T(g)(t)|^s dt
&\le n^s\int_0^1\|T(f(\tau))-T(g(\tau))\|_2^s d\tau\cr &\le
n^sA^s\int_0^1|f(\tau)-g(\tau)|^r d\tau,\end{array}$$ and
$$\int_0^1|\tilde T(f)(t)-\tilde T(g)(t)|^s dt\ge
\frac{1}{n^sA^s}\int_0^1|f(\tau)-g(\tau)|^r d\tau.$$ It follows that
$$\frac{1}{nA}\|f-g\|_r^\alpha\le\|\tilde T(f)-\tilde T(g)\|_s\le
nA\|f-g\|_r^\alpha.$$ Thus $\tilde T$ is a H\"older$(\alpha)$
embedding. Since $L_s\hookrightarrow\ell_q(L_s)$ and
$\alpha=\frac{p}{q}$, $L_r$ H\"older$(\frac{p}{q})$ embeds into
$\ell_q(L_s)$. Then $E(L_r,p)\le_B E(L_s,q)$ follows from Corollary
\ref{reduce1}.

Similarly, we can prove that $E(\ell_r,p)\le_B E(\ell_s,q)$.
\hfill$\Box$

\begin{theorem}\label{unreducible}
For $r,s,p,q\in[1,+\infty)$, if $E(\ell_r,p)\le_B E(\ell_s,q)$ or
$E(L_r,p)\le_B E(L_s,q)$, then we have
\begin{enumerate}
\item[(1)] $p\le q$;
\item[(2)]
$\min\left\{\frac{r}{p},\frac{2}{p}\right\}\ge\min\left\{\frac{s}{q},1,\frac{2}{q}\right\}$;
\item[(3)] $\max\{r,2\}\le\max\{s,q,2\}$.
\end{enumerate}
\end{theorem}

{\bf Proof.} Recall that
$$p(\ell_s)=p(L_s)=\min\{s,2\},\quad q(\ell_s)=q(L_s)=\max\{s,2\},$$
and for any Banach space $X$,
$$p(\ell_q(X))=\min\{p(X),q\},\quad q(\ell_q(X))=\max\{q(X),q\}.$$ Thus
$$p(\ell_q(\ell_s))=p(\ell_q(L_s))=\min\{s,q,2\},$$
$$q(\ell_q(\ell_s))=q(\ell_q(L_s))=\max\{s,q,2\}.$$

Therefore, clauses (1),(2) and the case $\min\{s,q\}>1$ in clause
(3) follow from Theorems~\ref{freduce},\ref{mtype},\ref{mcotype} and
Lemma~\ref{alpha}.

For the case $\min\{s,q\}=1$, let $\alpha\in(0,1)$ be arbitrary. By
Theorem~\ref{reducible}, (5), we have
$$E(\ell_s,q)\le_B E\left(\ell_\frac{s}{\alpha},\frac{q}{\alpha}\right),\quad
E(L_s,q)\le_B E\left(L_\frac{s}{\alpha},\frac{q}{\alpha}\right).$$
Thus
$\max\{r,2\}\le\max\left\{\frac{s}{\alpha},\frac{q}{\alpha},2\right\}$
for any $\alpha\in(0,1)$. It follows that
$\max\{r,2\}\le\max\{s,q,2\}$. \hfill$\Box$

\begin{corollary}
For $r,s,p,q\in[1,+\infty)$, if $E(\ell_r,p)\sim_B E(\ell_s,q)$ or
$E(L_r,p)\sim_B E(L_s,q)$, then we have $p=q$ and
$$r=s\quad\mbox{ or }\quad p\le r,s\le 2\quad\mbox{ or }\quad 2\le r,s\le p.$$
\end{corollary}

\begin{corollary}\label{iff}
For $r,s\in[1,2]$ and $p,q\in[1,+\infty)$, if $s\le q$, then
$$E(L_r,p)\le_B E(L_s,q)\iff p\le
q,\;\frac{r}{p}\ge\frac{s}{q}.$$
\end{corollary}

{\bf Proof.} ``$\Rightarrow$'' follows from
Theorem~\ref{unreducible}.

($\Leftarrow$). Assume that $2^k\le\frac{qr}{p}\le 2^{k+1}$ for some
$k\in\Bbb N$. From Theorem~\ref{reducible}, (4) and (5), we have
$$\begin{array}{ll}E(L_r,p)&\le_B E\left(L_2,\frac{2p}{r}\right)\le_B
E\left(L_1,\frac{2p}{r}\right)\le_B
E\left(L_2,\frac{4p}{r}\right)\le_B\cdots\cr &\le_B
E\left(L_2,\frac{2^jp}{r}\right)\le_B
E\left(L_1,\frac{2^jp}{r}\right)\le_B
E\left(L_2,\frac{2^{j+1}p}{r}\right)\le_B\cdots\cr &\le_B
E\left(L_2,\frac{2^kp}{r}\right).\end{array}$$ If
$s\ge\frac{qr}{2^kp}$, denote $s_1=\frac{2^kps}{qr}$. Then
$s_1\in[1,2]$, we have
$$E\left(L_2,\frac{2^kp}{r}\right)\le_B
E\left(L_{s_1},\frac{2^kp}{r}\right)\le_B E(L_s,q).$$ Otherwise,
denote $s_2=\frac{qr}{2^kp}$. Then $s_2\in[1,2]$ and $s<s_2$, we
have
$$E\left(L_2,\frac{2^kp}{r}\right)\le_B
E\left(L_1,\frac{2^kp}{r}\right)\le_B E(L_{s_2},q)\le_B E(L_s,q).$$
\hfill$\Box$

\begin{corollary}
For $r,s,p,q\in[1,+\infty)$, if $s\le q\le 2$, then
$$E(L_r,p)\le_B E(L_s,q)\iff p\le
q,\;r\le 2,\;\frac{r}{p}\ge\frac{s}{q}.$$
\end{corollary}

{\bf Proof.} ``$\Rightarrow$'' follows from
Theorem~\ref{unreducible}, and ``$\Leftarrow$'' follows from
Corollary~\ref{iff}. \hfill$\Box$

In the end, we settle down the case $X=c_0$ or $C[0,1]$.

\begin{theorem}
For $r,p,q\in[1,+\infty)$, we have
\begin{enumerate}
\item[(1)] $E(C[0,1],p)\sim_B E(c_0,p)$;
\item[(2)] $E(L_r,p)\le_B E(c_0,p)$;
\item[(3)] $q\in[1,+\infty)$, $E(c_0,p)\not\le_B E(L_r,q)$;
\item[(4)] if $p<q$, then $E(c_0,p)<_B E(c_0,q)$.
\end{enumerate}
\end{theorem}

{\bf Proof.} $E(c_0,p)\le E(C[0,1],p)$ is trivial. From Theorem
\ref{aharoni}, $C[0,1]\hookrightarrow c_0$ and $L_r\hookrightarrow
c_0$. So clauses (1) and (2) hold.

(3) Fix an $s>\max\{r,q,2\}$. Theorem~\ref{unreducible}, (3) shows
$E(L_s,p)\not\le_B E(L_r,q)$.  So $E(c_0,p)\not\le_B E(L_r,q)$,
since $E(L_s,p)\le_B E(c_0,p)$.

(4) From Lemma~\ref{fractal}, there are $n\in\Bbb N$ and a
H\"older$(\frac{p}{q})$ embedding $T:\Bbb R\to\Bbb R^n$. Fix a
bijection $\langle\cdot,\cdot\rangle_n:\Bbb
N\times\{1,2,\cdots,n\}\to\Bbb N$. We define $\hat T:c_0\to c_0$ by
$$\hat T(x)(\langle k,m\rangle_n)=T(x(k))(m)$$ for $k\in\Bbb N$ and
$m=1,2,\cdots,n$. It is easy to check that $\hat T$ is a
H\"older$(\frac{p}{q})$ embedding. It follows that $E(c_0,p)\le_B
E(c_0,q)$.

On the other hand, $E(c_0,q)\not\le_B E(c_0,p)$ follows from
Lemma~\ref{alpha}. \hfill$\Box$

\section{Further remarks}

Perhaps the most curious problem is how to compare equivalence
relations $E(X,p)$ and $E(Y,p)$ when $X$ and $Y$ have same types and
cotypes. Especially, if $r\ne 2$, does $E(L_r,p)\sim_B E(\ell_r,p)$?
Though $L_r\not\hookrightarrow\ell_r$, we have the following lemma.

\begin{lemma}
For $r\in[1,+\infty)$, $L_r$ finitely Lipschitz embeds into
$\ell_r$.
\end{lemma}

{\bf Proof.} Let $F=\{f_1,\cdots,f_n\}\subseteq L_r$. Denote
$\varepsilon=\min\{\|f_i-f_j\|_r:1\le i<j\le n\}>0$. Find continuous
functions $\varphi_1,\cdots,\varphi_n$ such that for each $i$, we
have $\|f_i-\varphi_i\|_r<\frac{\varepsilon}{8}$. Since all
$\varphi_i(t)$'s are uniformly continuous on $[0,1]$, there exists a
sufficiently large $m$ such that, for $k<m,i=1,\cdots,n$ and
$t\in\left[\frac{k}{m},\frac{k+1}{m}\right]$, we have
$$\left|\varphi_i(t)-\varphi_i\left(\frac{k}{m}\right)\right|<\frac{\varepsilon}{8}.$$
Thus for $i,j=1,\cdots,n$, we have
$$\left|\|\varphi_i-\varphi_j\|_r-\left(\sum_{k=0}^{m-1}\frac{1}{m}\left|\varphi_i\left(\frac{k}{m}\right)
-\varphi_j\left(\frac{k}{m}\right)\right|^r\right)^{\frac{1}{r}}\right|<\frac{\varepsilon}{4}.$$
Define $T_F:F\to\ell_r$ by
$$T_F(f_i)=m^{-\frac{1}{r}}\left(\varphi_i(0),\varphi_i\left(\frac{1}{m}\right),\cdots,\varphi_i\left(\frac{m-1}{m}\right),0,0,\cdots\right)$$
for $i=1,\cdots,n$. Then we can check that, for $1\le i<j\le n$,
$$\frac{1}{2}\|f_i-f_j\|_r\le\|f_i-f_j\|_r-\frac{\varepsilon}{2}\le\|T_F(f_i)-T_F(f_j)\|_r\le\|f_i-f_j\|_r+\frac{\varepsilon}{2}
\le2\|f_i-f_j\|_r,$$ as desired. \hfill$\Box$

\begin{question}
For $r,p\in[1,+\infty)$, if $r\ne 2$, does $E(L_r,p)\le_B
E(\ell_r,p)$?
\end{question}

If the statement in this question is true, it will follow that, for
$1\le p\le r\le 2$,
$$E(\ell_r,p)\le_B E(L_r,p)\le_B E(L_p,p)\le_B E(\ell_p,p)\le_B
E(\ell_r,p).$$ Then we shall have
$$\Bbb R^\Bbb N/\ell_p\sim_B E(\ell_p,p)\sim_B E(\ell_r,p)\sim_B E(L_r,p)\sim_B
E(L_2,p).$$

Though Corollary \ref{iff} gives an almost complete picture on Borel
reducibility between $E(L_r,p)$'s for $r\in[1,2]$, we know little
about the case $r\ge 2$. So another problem is, whether clauses
(1)-(3) in Theorem \ref{unreducible} can be a sufficient condition
for $E(L_r,p)\le_B E(L_s,q)$. This problem leads to the following
question.

\begin{question}
\begin{enumerate}
\item[(1)] For $r,s\ge 2$, if $r\le s$, does $L_r$ finitely Lipschitz embed
into $L_s$? Furthermore, does $E(L_r,p)\le_B E(L_s,p)$?
\item[(2)] For $p,q\in[1,+\infty),r\ge 2$, if $p\le q$, does $L_r$ finitely
H\"older$(\frac{p}{q})$ embed into $L_r$ itself? Furthermore, does
$E(L_r,p)\le_B E(L_r,q)$?
\end{enumerate}
\end{question}

\section{Acknowledgments}

I would like to thank Su Gao for leading me into the field of
descriptive set theory. I am grateful to Zhi Yin and Minggang Yu for
conversations in seminars. Special thanks are due to Guanggui Ding,
Guimei An and Lei Li for their help with the references in
functional analysis.

\label{}




\begin{thebibliography}{00}

\bibitem{aharoni} I. Aharoni, Every separable metric space is
Lipschitz equivalent to a subset of $c_0^+$, Israel J. Math. 19
(1974) 284-291.

\bibitem{BK} H. Becker, A. S. Kechric, The Descriptive Set Theory of
Polish Group Actions, London Math. Soc. Lecture Notes Series, vol.
232, Cambridge University Press, 1996.

\bibitem{BL} Y. Benyamini, J. Lindenstrauss, Geometric Nonlinear
Functional Analysis, Amer. Math. Soc. Colloq. Publ., vol. 48, A. M.
S., 2000.

\bibitem{BMW} J. Bougain, V. Milman, H. Wolfson, On type of metric
spaces, Trans. Amer. Math. Soc. 294 (1986) 295-317.

\bibitem{DH} R. Dougherty, G. Hjorth, Reducibility and
nonreducibility between $\ell^p$ equivalence relations, Trans. Amer.
Math. Soc. 351 (1999) 1835-1844.

\bibitem{farah} I. Farah, Basis problem for turbulent actions II:
$c_0$-equalities, Proc. London Math. Soc. (3) 82 (2001) 1-30.

\bibitem{gao} S. Gao, Equivalence relations and classical Banach
spaces, in: S.S. Goncharov, R. Downey, H. Ono (eds.), Mathematical
Logic in Asia, Proceedings of the 9th Asian Logic Conference,
Novosibirsk, Russia, 2005, World Scientific, 2006, pp. 70-89.

\bibitem{gaobook} S. Gao, Invariant Descriptive Set Theory,
Monographs and Textbooks in Pure and Applied Mathematics, vol. 293,
CRC Press, 2008.

\bibitem{heinrich} S. Heinrich, Ultraproducts in Banach space
theory, J. Reine Angew. Math. 313 (1980) 72-104.

\bibitem{hjorth} G. Hjorth, Actions by the classical Banach spaces,
J. Symb. Logic 65 (2000) 392-420.

\bibitem{kanovei} V. Kanovei, Borel Equivalence Relations: Structure
and Classification, University Lecture Series, vol. 44, A. M. S.,
2008.

\bibitem{LT} J. Lindenstrauss, L. Tzafriri, Classical Banach Spaces,
II: Function Spaces, Ergebnisse der Mathematik und ihrer
Grenzgebiete, vol. 97, Springer-Verlag, 1979.

\bibitem{MP} B. Maurey, G. Pisier, S\'eries de variables al\'eatoires
vectorielles ind\'ependantes et propri\'et\'es g\'eom\'etriques des
espaces de Banach, Studia Math. 58 (1976) 45-90.

\bibitem{MN} M. Mendel, A. Naor, Metric cotype, Ann. Math. 168
(2008) 247-298.

\bibitem{pisier} G. Pisier, Holomorphic semigroups and the geometry
of Banach spaces, Ann. Math. 115 (1982) 375-392.

\bibitem{T-J} N. Tomczak-Jaegermann, Banach-Mazur Distances and
Finite-Dimensional Operator Ideals, Pitman Monographs and Surveys in
Pure and Applied Mathematics, vol. 38, Longman Scientific \&
Technical, 1989.

\end{thebibliography}
\end{document}